\documentclass{article}
\usepackage{amsmath}
\usepackage{amsthm}
\usepackage{latexsym}
\usepackage{fullpage}
\usepackage{graphicx}
\usepackage{epstopdf}
\usepackage{amsfonts}
\usepackage{mathtools}
\usepackage{wrapfig}
\usepackage{stix}
\usepackage{url}
\usepackage{hyperref}
\usepackage{amsrefs}
\usepackage[english]{babel}  
\usepackage{wrapfig}

\newcommand{\Frac}[2]{\displaystyle\frac{#1}{#2}}

\newtheorem{definition*}{Definition:}
\newtheorem{facts*}{Facts:}

\begin{document}
\setlength{\parindent}{.5cm}
\setlength{\parskip}{.2cm}
{\title{A model for the space of convex quadrilaterals.}}
{\author{Carl Eberhart \\ \textit{carl@g.uky.edu}}}
\maketitle
\begin{abstract}
\noindent This paper describes a 'model' for the space of convex quadrilaterals, that is, a set $\mathbb{Q}$ of convex quadrilaterals in the plane which gives a 'cross-section' of the equivalence classes of similar convex quadrilaterals, in the sense  that every convex quadrilateral is similar to exactly one member of $\mathbb{Q}$.  It is also has the property that quadrilaterals which have nearly the same vertices are close to each other in the model in the Hausdorff metric.

\noindent We also parameterize our model with the 4-cell and use that to describe and investigate various classes of quadrilaterals.  For example, we show that the trapezoids form a 3-dimensional closed subset of the model which separates the model into 3 disjoint open  sets, each homeomorphic with a 4-cell.  We used \textbf{SageMath} to compose the latex for this note, and to make the \textbf{Sage Cell Interacts} to explore the model, available at \textit{sagelets.org}.

 Quadrisections are discussed and a question is posed about the maximum finite number of quadrisections of a convex quadrilateral.

\end{abstract}

\section{Introduction}

  Recall that a \textit{convex quadrilateral} is a polygon with 4 vertices such that the intersection of the diagonals lies in the interior of each .  There are many different kinds of convex quadrilaterals: \textit{trapezoid}, \textit{parallelogram}, \textit{rhombus}, \textit{rectangle}, \textit{kite}, \textit{cyclic quadrilateral} and many others. Over the centuries, there have been numerous classification schemes proposed.  Martin Josefsson~\cite{Josefsson} has given a good summary of these in his paper  \textit{On the classification of convex quadrilaterals} and has proposed another interesting classification.
In addition,  there have been studies done of various spaces of  polygons and quadrilaterals in particular.

In this paper we describe a 'model' for the space of convex quadrilaterals, that is, a set $\mathbb{Q}$  of convex quadrilaterals in the
plane which gives a 'cross section' of the equivalence classes of similar convex quadrilaterals in the sense  that every convex quadrilateral is similar to exactly one member of $\mathbb{Q}$.  This model also has the property that quadrilaterals which have nearly the same vertices are close to each other.  However, it (unavoidably) also has nearly congruent quadrilaterals which are far apart.

  We modeled the space of triangles in \cite{Eberhart}. By scaling and a Euclidean motion, we can place each triangle so that the its vertices are $A=(1,0)$, $B=(x,y)$, and $C=(-1,0)$, where $x\ge 0$, $y>0$, and $(x+1)^2+y^2 \le 4$, that is, $B$ is in first quadrant above the $x$-axis and  $|BC|\le 2$, where $|PQ|$ is the Euclidean distance between points $P$ and $Q$.  That model is 2 dimensional and is  a disk with a closed arc removed from the boundary  (see the figure \textbf{Space of triangles}).  We also investigated how the various types of triangles were situated in the model.

\begin{wrapfigure}{R}{0.6\textwidth}
\centering
\includegraphics[scale=.4]{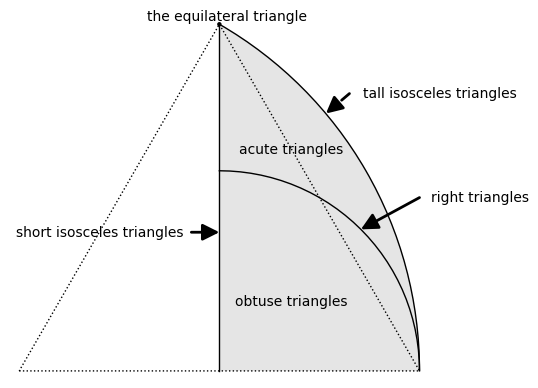}\\\textbf{Space of triangles.}
\end{wrapfigure}
  For example, the boundary of the model consists  of an open interval of isosceles triangles with the equilateral triangle at the midpoint with the tall isosceles triangles to the right and the short isosceles triangles to the left.  Also, the right triangles forms an arc (half-open) ending at the isosceles right triangle and separating the model into two components, one being all triangles with an obtuse angle and the other being all triangles with no obtuse angle.

  We were primarily interested in identifying the triangles with one, two or three \textit{quadrisections} (perpendicular segments dividing a triangle into four equal areas), and found that almost all triangles have only one quadrisection save a small open (in the model) set of triangles with three quadrisections  about the equilateral triangle.  The boundary of this  is a closed arc separating the triangle with three quadrisections from those with only one quadrisection, starting at the only isosceles triangle with two quadrisections and consisting of scalene triangles with two quadrisections except for the other endpoint which is an isosceles triangle with one quadrisection.

A little thought tells us that any model for the convex quadrilaterals is 4 dimensional:  jJust fix two vertices and let the other two roam subject to the conditions that the quadrilateral must be convex and not degenerate to a triangle or segment. Also we don't want to have two congruent quadrilaterals in our model, but we want every quadrilateral to be similar to one quadrilateral in our model.

\section{Description of the Model}\label{description}

We will build our model for the convex quadrilaterals from the set $\mathbb{Q}$ of quadrilaterals $ABCD$ which satisfy the following conditions:

\textbf{Conditions for membership in $\mathbb{Q}$}:
\textit{$A=(1,0)$, $C=(-1,0) $, $B=(x,y)$ with $x\ge 0, y>0$, and $D=(z,w)$ with $w <0$, $z\le x$ ,and  $(x-z)^2+(y-w)^2 \le 4$.  Also, we require that the intersection $X=(u,0)$ of $AC$ and $BD$ satisfy $0\le u < 1$.}\label{conds1}

\textit{In other words, for each member $ABCD$ of $\mathbb{Q}$, the diagonal $BD$ has length no more than 2, $AX$ has length no more than 1, and $\angle AXB$ is no more than 90 degrees.}

If $ABCD$ is in $\mathbb{Q}$, we call $AC$ and $BD$ the \textbf{major} and \textbf{minor} diagonals of $ABCD$.

We claim that each convex quadrilateral Q is similar to one or possibly two members of $\mathbb{Q}$.   By scaling, rotating, reflecting and translating we can move the quadrilateral so that one diagonal is $AC$ with $A=(1,0)$ and $C=(-1,0)$ in eight, four, two, or one orientations, depending on its symmetries. One of these orientations, or rarely two of them, satisfy the conditions for membership in $\mathbb{Q}$ above.

In the rare cases where there are duplicates copies of Q in $\mathbb{Q}$, we establish rules for which one to remove from $\mathbb{Q}$ below.

Here are typical examples of quadrilaterals with eight or four orientations.

\begin{table}[ht]
\centering
\begin{tabular}{cc}
\includegraphics[scale=.4]{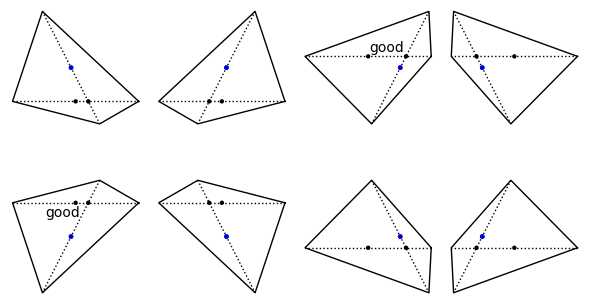}&\includegraphics[scale=.5]{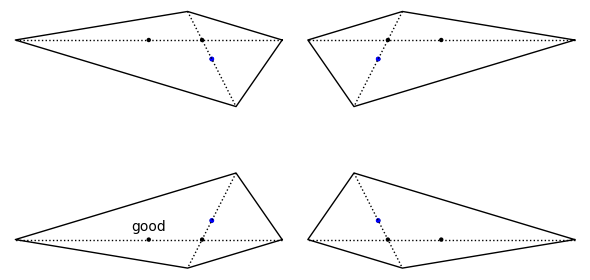}\\
\textbf{Eight orientations}&\textbf{Four orientations}
\end{tabular}
\end{table}

The duplicates occur when (1) the diagonals are perpendicular,  (2) the minor diagonal bisects the major diagonal, or (3) the diagonals have the same length.  Each case occurs on the \textbf{boundary} of $\mathbb{Q}$ and there is only one duplicate and there are good rules for elimination of one of the duplicates.  A member $ABCD$ is in the \textbf{boundary} of $\mathbb{Q}$ if it is possible to move $ABCD$ out of $\mathbb{Q}$ by moving $B$ or $D$ by an arbitrarily small distance. The vast bulk of the quadrilaterals in $\mathbb{Q}$ lie in the \textbf{interior} of $\mathbb{Q}$ consisting of those $ABCD$'s in $\mathbb{Q}$ for which there is a positive number $\epsilon$ such that each quadrilateral $AB'CD'$ with $A=(1,0)$, $C=(-1,0)$, and $B'$, $D'$ within $\epsilon$ of $B$, $D$ respectively lies in $\mathbb{Q}$.

The boundary of $\mathbb{Q}$ consists of three special classes of quadrilaterals.  (1) \textbf{orthodiagonal} quadrilaterals, those where the diagonals are perpendicular, (2) \textbf{bimajor} quadrilaterals, those whose minor diagonal \textbf{bi}sects the \textbf{major} diagonal, and (3) \textbf{equidiagonal} quadrilaterals, those with diagonals of the same length.   The first and last classes are named in Josefsson's\cite{Josefsson} classification scheme.  The term bimajor is of our making.

There is another class which has not been named before which is important in our discussion, namely, the class of quadrilaterals whose minor diagonal is bisected by its major diagonal.  This class contains the parallelograms and the kites, but is much larger than either.  We agree to call this class the \textbf{biminor} quadrilaterals\label{biminor}.

If $ABCD \in \mathbb{Q}$ is orthodiagonal, then it has a \textbf{twin} $AB'CD'$ obtained by reflecting $ABCD$ about the $x$-axis.
These are congruent orthodiagonl quadrilaterals in $\mathbb{Q}$, and one must be deleted.
  We could eliminate either of these from $\mathbb{Q}$, but must make a consistent choice.  Our rule below will eliminate the one which sticks up further above the $x$-axis, that is, the one where $|BX| > |DX|$. Note that if $ABCD$ is it's own twin (ie $|BX| = |DX|$), then it is a \textbf{kite}.

\begin{center}
\includegraphics[scale=.4]{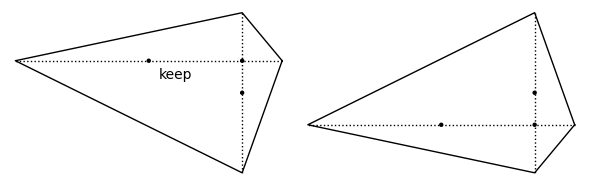}\\\textbf{Duplicate orthodiagonal quadrilaterals in $\mathbb{Q}$}
\end{center}

For bimajor quadrilaterals there are pairs $ABCD$ and $AB'CD'$ which are congruent under a rotation of 180 degrees about the origin.  Here again we will eliminate the quadrilateral where $|BX| > |DX|$. If $|BX|=|DX|$  (ie, $ABCD$ is biminor), then $ABCD$ is a \textbf{parallelogram}.

\begin{center}
\includegraphics[scale=.4]{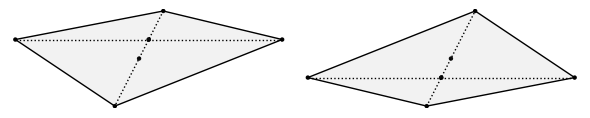}\\\textbf{Duplicate bimajor quadrilaterals in $\mathbb{Q}$}
\end{center}

Certain equidiagonal quadrilaterals in $\mathbb{Q}$ also have duplicates in $\mathbb{Q}$.  There are of two types.  (1) If $|AX|<|BX| < |DX|$. Then by reflecting $ABCD$ about the line which bisects $\angle BXA$ and following with a horizontal translation we can move $BD$ onto $AC$ yielding $AB'CD' \in \mathbb{Q}$ congruent to $ABCD$ with $|AX'|>|B'X'|$.  We will eliminate from $\mathbb{Q}$ the one with $|AX|>|BX|>$.  If $|BX|=|AX|$, then $|DX|=|CX|$ and by side angle side, $ABCD$ is an isosceles trapezoid, what we term later a \textbf{tall} isosceles trapezoid.  (2)  $|BX| > |DX|$. In this case, reflection through the bisector of $\angle BXA$ puts $B'$ onto the wrong side of $X$, but rotation of 180 degrees about $(0,0)$ followed by a horizontal translation carries $B'D'$ onto $CA$. Here we agree to eliminate the quadrilateral with $|DX|>|AX|$
Again, if $|AX|$ and $|DX|$ have the same length, there is only one member of $\mathbb{Q}$ and it is a  (\textbf{short} as we term later) \textbf{isosceles quadrilateral}.

\textbf{Note:} There is a sagelet in \textit{sagelets.org} to check the duplications.

\begin{center}
\begin{tabular}{cc}
\includegraphics[scale=.4]{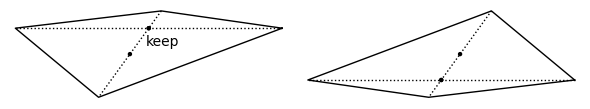} & \includegraphics[scale=.4]{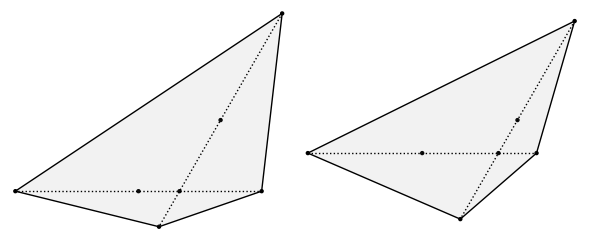}\\
\textbf{Equidiagonal duplicates (1)} & \textbf{Equidiagonal duplicates (2) }
\end{tabular}
\end{center}

To summarize, we have these

\textbf{Elimination Rules:} (1) If $ABCD \in \mathbb{Q}$ is orthodiagonal and $|BX|> |DX|$, then eliminate $ABCD$ from $\mathbb{Q}$.   (2) If it is bimajor and $|BX|> |DX|$, then eliminate $ABCD$ from $\mathbb{Q}$.  (3) If it is equidiagonal with  $|AX| < |BX| < |DX|$ or $|AX|< |DX|< |BX|$, then eliminate $ABCD$ from $\mathbb{Q}$.

At this point,  the description of the model $\mathbb{Q}$ is complete.   Now we will \textbf{parameterize} it and use the \textbf{parameter space} to help visualize the various classes of quadrilaterals and their relation to each other.

\section{Parameterizing  the Model.}

\begin{wrapfigure}{R}{0.6\textwidth}
\centering
\includegraphics[scale=.5]{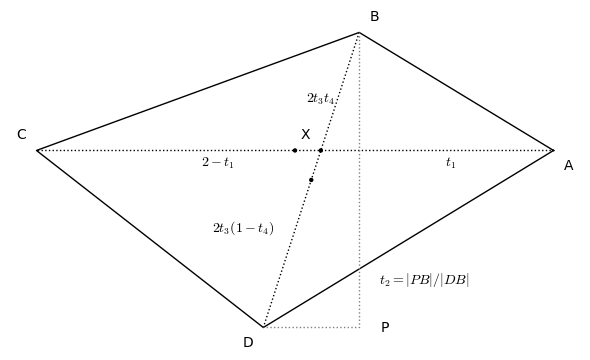}\\\textbf{Parameters determine $ABCD$.}
\end{wrapfigure}
From \ref{description}, recall that a typical element of $\mathbb{Q}$ has vertices $A=(1,0)$, $B=(x,y)$, $C=(-1,0)$, and $D=(z,w)$ where $AC$ and $BD$ intersect at $X=(u,0)$ and $x,y,z,w,u$ satisfy the restrictions $(x-z)^2+(y-w)^2 \le 4$, $0,z \le x$, $w< 0< y$,  and $0 \le u < 1$. In addition, by the elimination rules above, if $|BD|=2$ then $|BX| \le |AX|$, and if $\angle AXB$ is a right angle, then $|BX| \le |DX|$.

We could use the coordinates of $B$ and $D$, $x,y,z,w$ to parameterize $ABCD$, but these are very inconvenient parameters to work with.  We prefer working with parameters in the unit interval $[0,1]$.

The most convenient set of parameters for a member $\mathcal{Q}=ABCD$ of $\mathbb{Q}$ we have found  comes from the major and minor diagonals $AC$ and $BD$ of $\mathcal{Q}$ and their intersection point $X$.  Associate with $\mathcal{Q}$ a unique 4-tuple $(t_1,t_2,t_3, t_4)$ with $t_i \in (0,1], i\in \{1,2,3,4\}$ as follows:

 $\mathbf{t_1 =  \Frac{2\,|AX|}{|AC|}}$.  So $t_1 \in (0,1]$ when $X$ is in the right hand half of $AC$. 

  $\mathbf{t_2=\Frac{\angle BXA}{90}}$,  $\angle BXA$ is measured in degrees.  So $t_2 \in (0,1]$ determines $\angle AXB$.

  $\mathbf{t_3=\frac{1}{2}|BD|}$.  So since $2\,t_3=|BD|\le |AC| =2$, $t_3 \in (0,1]$.

 $\mathbf{t_4=\Frac{|BX|}{|BD|}}$, so $t_4 \in (0,1)$. Note that $|DX| = 2\,t_3\,(1-t_4)$.

Note that $t_1=1$ if and only if $\mathcal{Q}$ is bimajor,  $t_2=1$ if and only if  $\mathcal{Q}$ is orthodiagonal, $t_3=1$ if and only if  $\mathcal{Q}$ is equidiagonal and $t_4=.5$ if and only if  $\mathcal{Q}$ is biminor.

This sets up a 1-1 correspondence $f:\mathcal{Q}=ABCD \to (t_1,t_2,t_3,t_4)$  from $\mathbb{Q}$ to the \textbf{parameter space}  $\mathbb{P}$ of $\mathbb{Q}$, contained in the  4-cell $\mathbb{I}^4=(0,1]\times (0,1]\times (0,1]\times (0,1)$. The correspondence is not onto $\mathbb{I}^4$, because of the elimination of certain duplicate equidiagonal, orthodiagonal, and bimajor quadrilaterals from $\mathbb{Q}$.

We eliminated the bimajor quadrilaterals from $\mathbb{Q}$  with $2\,t_3\,t_4=|BX|>|AX|=1$ or $t_4 \ge t_3\,t_4 > \Frac{1}{2}$, so points $ (1,t_2,t_3,t_4)\in \mathbb{I}^4$  with $  \Frac{1}{2}< t_4 $  are not in $\mathbb{P}$.
Also, we have eliminated the orthodiagonal quadrilaterals $ABCD$ with  $2\,t_3\,t_4=|BX|> |DX|=2\,t_3\,(1-t_4)$. So points $ (t_1,1,t_3,t_4)\in \mathbb{I}^4$  with $t_4<\Frac{1}{2}$ are not in $\mathbb{P}$.
Finally, we eliminated the duplicate equidiagonal quadrilaterals ( $t_3=1$) in $\mathbb{Q}$ with $t_1=|AX| < |BX|= 2\,t_4 < |DX|=2(1-t_4)$ or $|AX|=< |DX|< |BX|$ .  So points $(t_1,t_2,1,t_4)\in \mathbb{I}^4$ with  $ \Frac{t_1}{2}< t_4 < 1-\Frac{t_1}{2}$ are not in $\mathbb{P}$.

So the parameter space $\mathbb{P}$ of $\mathbb{Q}$ is $ \mathbb{I}^4 - (\mathbb{B}_1 \cup \mathbb{B}_2 \cup \mathbb{B}_3)$, where  $\mathbb{B}_1$ is
one-half of the 1-boundary in the $t_1$ coordinate, $\mathcal{B}_1= \{1\} \times (0,1] \times (0,1] \times (0,.5)$,
$\mathbb{B}_2$ is one half the 1-boundary in the $t_2$ coordinate, $\mathcal{B}_2= (0,1] \times \{1\} \times (0,1] \times (0,.5)$,
and  $\mathbb{B}_3$ is a strange looking one half  of the 1-boundary in the $t_3$ coordinate,
$\mathbb{B}_3= \bigcup_{t_1 \in (0,1]} \{t_1\} \times (0,1]\times \{1\} \times (0,t_1/2)$.  See the diagram \textbf{Graph of the boundary of $\mathbb{P}$}.

So now the correspondence $f:\mathbb{Q} \to \mathbb{P}$ is 1-1 \textit{and} onto.  It is also true that $f$ is a homeomorphism, that is, $f$ and $f^{-1}$ are continuous, where we are endowing the set $\mathbb{Q}$ with the \textbf{Hausdorff distance} (see Wikipedia) which simplifies in this case to the distance between  $ABCD$ and $AB'CD'$ in $\mathbb{Q}$ is the \textbf{maximum} of the distances $|BB'|$ and $|DD'|$.

Now quadrilateral can be thought of as a point in the 4-cell. The various classes of quadrilaterals and their relations to one another can be visualized by finding equations and/or inequalities in the parameters which characterize them.  These  algebraic characterizations enable us to determine their dimension and whether they separate $\mathbb{P}$ into disjoint pieces. We can also draw some pictures.

\textbf{Notations for the \textbf{graph} of the algebraic description  of  a class of quadrilaterals in $\mathbb{Q}$:}
\textit{If $\mathcal{C}$ is a set of equations and/or inequalities in the parameters $t_i, i\in\{1,2,3,4\}$ then  the points in $\mathbb{P}$ which satisfy them is  $\mathbb{P}(\mathcal{C})$}, and $\mathbb{Q}(\mathcal{C})$ is the class of quadrilaterals $\mathbf{Q} \in \mathbb{Q}$ such that $f(\mathbf{Q})\in \mathbb{P}(\mathcal{C})$.

So for example, $\mathbb{P}(t_2=1)$ is the graph of orthodiagonal quadrilaterals, $\mathbb{P}(t_3=1)$ is the graph of equidiagonal quadrilaterals, and $\mathbb{P}(t_2=1)\cap\mathbb{P}(t_3=1)=\mathbb{P}(t_2=1,t_3=1)$ is called the \textbf{midsquare quadrilaterals} in Josefsson's terminology\cite{Josefsson},p.81. The class of the graph $\mathbb{P}(t_1=1)$ is the bimajor quadrilaterals.   In any case, we see that \textit{any class whose graph is contained in one or more of $\mathbb{P}(t_i=1)$, $i \in \{1,2,3\}$ lies in the boundary of $\mathbb{P}$}.
\begin{center}
\includegraphics[scale=.4]{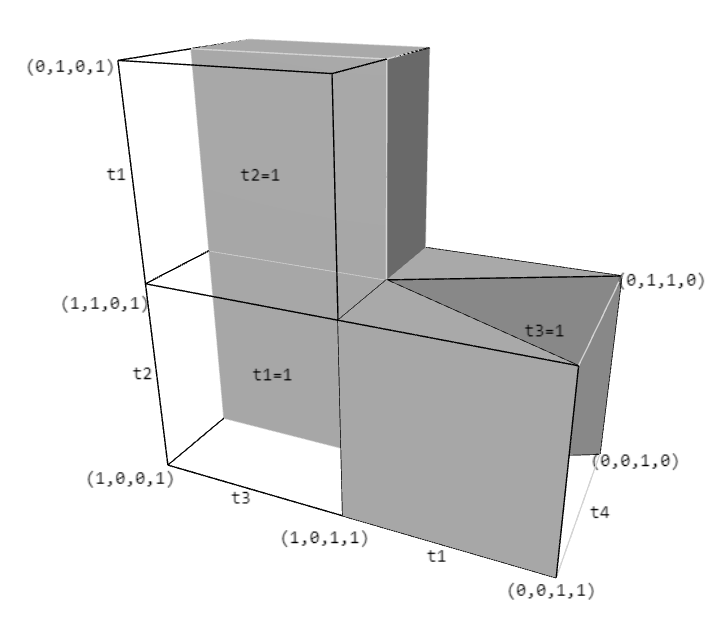}\\\textbf{Graph of the boundary of  $\mathbb{P}$}
\end{center}

\subsection{Trapezoids}

There are two  types of trapezoids in $\mathbb{Q}$, depending on whether $AB$ and $CD$  or $BC$ and $AD$ are parallel, These will be  referred to these as \textbf{tall} and \textbf{short}  trapezoids respectively.

You can see using alternate interior angles that $\triangle ABX$ and $\triangle CDX$ are similar for tall trapezoids and $\triangle BCX$ and $\triangle DAX$ are similar for short trapezoids. (see figure \textbf{Tall and Short trapezoids} )
\begin{center}
\includegraphics[scale=.6]{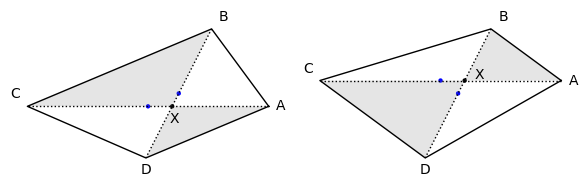}\\\textbf{Short and Tall trapezoids}
\end{center}

If $ABCD \in \mathbb{Q}$ is a tall trapezoid,  then $\Frac{|AX|}{|CX|}=\Frac{|BX|}{|DX|}$, or $\Frac{t_1}{2-t_1}=\Frac{2\,t_3\,t_4}{{2\,t_3\,(1-t_4)}}$.  This simplifies to $t_4=\Frac{t_1}{2}$, which gives us a description of the tall trapezoids in $\mathbb{Q}$ in terms of its parameters. Likewise, short trapezoids are characterized by the equation $t_4=1-\Frac{t_1}{2}$.

 The remaining quadrilaterals (which form the vast majority) constitute 3 disjoint topological open 4-cells whose closures
 we will call \textbf{short, middle} and \textbf{tall  quadrilaterals} respectively.
This means that \textit{in order to deform a short quadrilateral into a tall quadrilateral through an arc of quadrilaterals in $\mathbb{Q}$, it is necessary to pass through the middle quadrilaterals}.

We can also work out descriptions of the parallelograms, rectangles, rhombi, the square, equidiagonal quadrilaterals, and orthodiagonal quadrilaterals.  From that, the dimension of each class is found and a diagram or picture of the class can be drawn in the parameter space. See the table \textbf{Algebraic description and dimension of some classes of quadrilaterals} at the end of this section.

The graphs which are 1, 2 or 3 dimensional can be visualized.  For example the graph of the parallelograms, $\mathbb{P}(t_1=1,t_4=.5)$ is dimension 2, in fact it is the 2-cell $\{1\}\times (0,1] \times (0,1] \times \{.5\} \subset \mathbb{P}$. It sits on the equidiagonal boundary ($t_1=1$) of the  quadrilaterals, and is the intersection of the tall and short trapezoids, each of which is a 3-cell having the parallelograms as a common face.  The 3 or 4 dimensional graphs can be visualized by looking at \textbf{cross-sections} along any one of the 4 parameters.  The graphs of the equidiagonal, orthodiagonal and bimajor quadrilaterals are examples of cross-sections of $\mathbb{P}$.  Another example is the graph of the  quadrilaterals $\mathbb{P}(t_4=1/2)$.  This class \textbf{separates} $\mathbb{P}$ into two disjoint open pieces.

In the figure \textbf{The space of quadrilaterals }, each point represents a 2-cell of quadrilaterals, for example, the 2-cell of parallelograms sits above the point labeled parallelograms. The cross-section of biminor quadrilaterals bisects the middle quadrilaterals into two pieces.  In fact, the mapping
$$f:\mathbb{P}\setminus \mathbb{P}(t_1=1\text{ or }t_2=1\text{ or }t_3=1)  \to \mathbb{P}\setminus \mathbb{P}(t_1=1\text{ or }t_2=1\text{ or }t_3=1)$$  given by $f(t_1,t_2,t_3,t_4)=(t_1,t_2,t_3,1-t_4)$ is a reflection, that is, a homeomorphism which is its own inverse.  It maps tall quadrilaterals to short and short to tall, with fixed point set the biminor quadrilaterals.

\begin{center}
\begin{tabular}{cc}
\includegraphics[scale=.4]{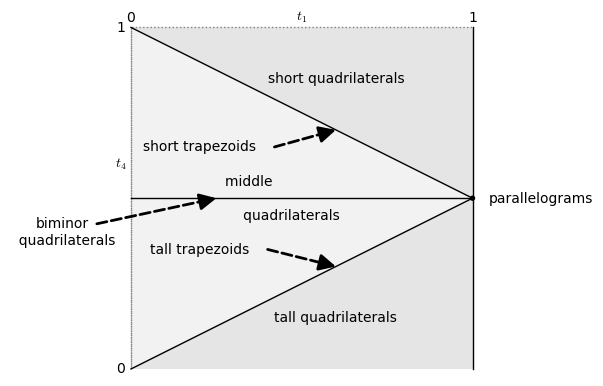} & \includegraphics[scale=.4]{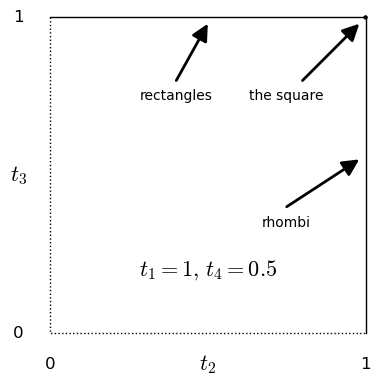}\\
\textbf{The space of quadrilaterals} & \textbf{Parallelograms }
\end{tabular}
\end{center}

The table below summarizes the classes of quadrilaterals, their algebraic characterizations, and their dimensions that we have discussed so far.

\begin{center}
\textbf{Algebraic description and dimension of some classes of quadrilaterals}

\begin{tabular}{|l|l|}
\hline
tall quadrilaterals: $t_4 \le \frac{1}{2}\,t_1$, dim 4  & short quadrilaterals: $t_4 \ge 1-\frac{1}{2}\,t_1$, dim 4  \\
    middle quadrilaterals $\frac{1}{2}\,t_1 \le t_4 \le 1-\frac{1}{2}\,t_1$, dim 4 &  interior quadrilaterals:  $t_1,t_2,t_3,t_4 \in (0,1)$ \\
  equidiagonal quadrilaterals: $t_3 = 1$,  dim 3 & orthodiagonal quadrilaterals: $t_2=1$, dim 3   \\

  bimajor quadrilaterals:  $t_1=1$,  dim 3 & biminor quadrilaterals:  $t_4= \frac{1}{2}$, dim 3 \\
 tall trapezoids:  $t_4 =  \frac{1}{2}\,t_1$,  dim 3 & short trapezoids: $t_4 = 1- \frac{1}{2}\,t_1$ , dim 3   \\

  kites:  $t_2=1$, $t_4=\frac{1}{2}$, dim 2 & midsquare quadrilaterals:  $t_2=1$, $t_3=1$, dim 2 \\
  tall isosceles quadrilaterals:  $t_3=1$, $t_4=\frac{1}{2}t_1$, dim 2 & short isosceles quadrilaterals:  $t_3=1$, $t_4=1-\frac{1}{2}t_1$, dim 2\\
parallelograms: $t_1=1$, $t_4=\frac{1}{2}$, dim 2 & rectangles:   $t_1,t_3=1$, $t_4=\frac{1}{2}$, dim 1     \\
  rhombi:   $t_1,t_2=1$,$t_4=\frac{1}{2}$, dim 1 &  the square:    $t_1,t_2,t_3=1$, $t_4=\frac{1}{2}$,dim 0      \\

\hline
\end{tabular}
\end{center}

Next, let's examine some classes whose descriptions in the parameters are more complicated.

\section{Circles and quadrilaterals.}
There are at least three ways a circle can interact with a quadrilateral $ABCD$:  (1) It might pass through $A,B,C,D$, or (2) it might be tangent to all four sides of $ABCD$ in which case its interior would be in the interior of $ABCD$, or (3) it might be tangent to all four lines $AB$,$BC$,$CD$, and $DA$ and in the exterior of $ABCD$.  This determines three classes of quadrilaterals: \textit{Cyclic, Tangential, and Extangential} quadrilaterals.\cite{Josefsson}.  There are algebraic descriptions in terms of $t_1, t_2, t_3$, and $t_4$, which will enable us to determine the dimension and separation properties of these and other classes of quadrilaterals.  It  turns out that in the tangential and extangential cases, the equations are not useful for finding specific quadrilaterals, and we introduce a new set of parameters to do that.

\subsection{Cyclic quadrilaterals}.  We know that all angles inscribed in a circle and subtended by the same chord are equal, since an inscribed angle of a circle is equal to one half of the central angle subtended by the chord.  As a consequence, if $ABCD \in \mathbb{Q}$ is cyclic, $\triangle AXB$ is similar with $\triangle DXC$.  So $\Frac{t_1}{2 t_3(1-t_4)}=\Frac{2 t_3t_4}{2-t_1}$. Solve this for $t_1$ to get $t_1=1-\sqrt{1-4 t_3^2t_4(1-t_4)}$, that is, $\mathbb{P}(t_1=1-\sqrt{1-4 t_3^2t_4(1-t_4)})$ is 3-dimensional.  So the cyclic quadrilaterals form a 3-dimensional class in $\mathbb{Q}$.

Do the cyclic quadrilaterals separate the quadrilaterals?Yes,
\textit{$\mathbb{Q}$ is separated into two disjoint open sets by the cyclic quadrilaterals.}  Here's a plausibility argument for this.  First, for each $ABCD$ in $\mathbb{Q}$, let $\text{Circum}(ABCD)$ be the unique circle which passes through the three vertices $A$,$B$, and $C$ of $ABCD$.  Decompose $\mathbb{Q}$
into three disjoint sets: $\mathbb{Cin},\mathbb{Cyc},\mathbb{Cout}$ where $ABCD\in\mathbb{Cin}$, $\mathbb{Cyc}$,or $\mathbb{Cout}$ according to whether the vertex $D$ lies inside, on, or outside $\text{Circum}(ABCD)$. So $\mathbb{Cyc}$ is the class of cyclic quadrilaterals, and it is clear that if a quadrilateral from $\mathbb{Cin}$ is deformed along an arc of quadrilaterals to one in $\mathcal{Cout}$, then it must pass through $\mathbb{Cyc}$.

We can visualize the graph of the cyclic quadrilaterals by drawing a typical $t_2=\text{const}$ cross section of the graph.
(See the drawing below). Then the graph is just the direct product of $(0,1]$ with the cross-section. It is clear that the cyclic quadrilaterals separate the middle quadrilaterals into two disjoint open sets.

\begin{center}
\begin{tabular}{cc}
\includegraphics[scale=.4]{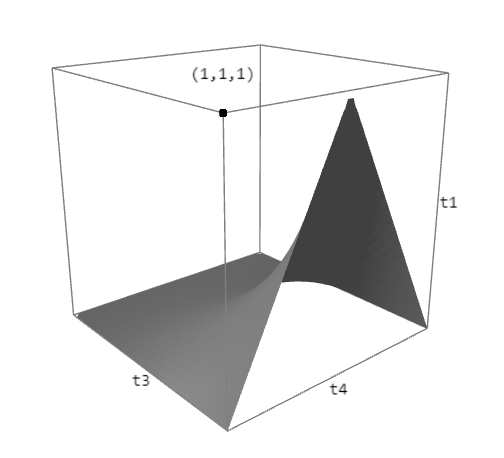} & \includegraphics[scale=.4]{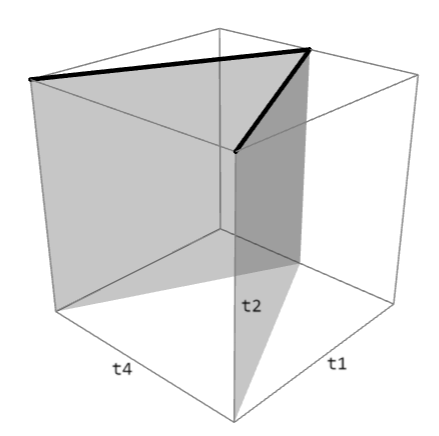}\\
\textbf{ $t_2 = \text{const}$ section of the cyclic quadrilaterals} & \textbf{ Graph of the isosceles trapezoids. }
\end{tabular}

\end{center}

What about the intersection of the trapezoids with the cyclic quadrilaterals?
Use a little algebra to show $\Frac{t_1}{2 t_3(1-t_4)}=\Frac{2 t_3t_4}{2-t_1}$ and ($t_1=2t_4$ or $t_1=1-2t_4$) implies $t_3=1$, so the intersection is the equidiagonal trapezoids, that is, the isosceles trapezoids.

Even though we can tell from the parameters of a quadrilateral if it is cyclic, we still don't know the center $O$ or radius $r$ of its circumcircle.  But since the perpendicular bisector of any chord passes through the center, we do know that $O=(0,h)$ for some $h$,
and from that $r=|O-A|=\sqrt{1+h^2}$.  But also $r=|O-B|$, yielding the equation $1+h^2 = x^2+(y-h)^2$ with $x=1-t_1+2t_3t_4 c_2$ and $y=2t_3t_4 s_2$, $s_2=\sin{(t_2\pi/2)}$, $c_2=\cos{(t_2\pi/2)}$.  Solve for $h$ to get $$h=\Frac{x^2+y^2-1}{2y}=\Frac{(1-t_1)^2+(2t_3t_4)^2+4(1-t_1)t_3t_4 c_2-1}{4t_3t_4 s_2}$$.

\subsection{Tangential and extangential quadrilaterals}

The tangential quadrilaterals, the ones with an \textbf{inscribed circle}, are characterized by the simple equation $|AB|+|CD|=|BC|+|AD|$ (sums of opposite sides are equal). You can verify this equation by dropping the center of the incircle perpendicularly onto each side, and noting that each side of the equation decomposes into rearrangements of the same four numbers.  The extangential quadrilaterals have a circle exterior to the quadrilateral which is tangent to all four extended sides. They were shown by Jacob Steiner in 1846 (see Wikipedia) to  be characterized by the equations ($|AB|+|BC|=|CD|+|AD|$ or $|AB|+|AD|=|CB|+|CD|$).  For extangential quadrilaterals in our model, this means that the exterior circle lies to the right of the major diagonal $AC$ if $|AB|+|BC|=|CD|+|AD|$ and above the minor diagonal $BD$ if $|AB|+|AD|=|CB|+|CD|$. We will term these \textbf{major extangential} and \textbf{minor extangential} quadrilaterals respectively.

An interesting observation about extangential quadrilaterals:  Each extangential quadrilateral $ABCD$ in $\mathbb{Q}$ determines an ellipse $\Frac{x^2}{(k/2)^2}+\Frac{y^2}{(k/2)^2-1}=1$ where $k=|AB|+|BC|$.  $A$ and $C$ are the foci of the ellipse and $B$ and $D$ lie on the ellipse. When $|BD|<2$, the quadrilateral is a major extangential quadrilateral.  Not all pairs $B$ and $D$ yield and extangential quadrilateral.  The $x$ coordinate of $B$ must be greater than the $x$ coordinate of $D$, and $B$ (and $D$) can't lie too far to the right (left):   the exact ranges have not been worked out.   If $|BD|>2$, then the diagram should be reflected about the bisector of $\angle BXA$ and rescaled and the foci of the ellipse become $B$ and $D$ and the quadrilateral changes to a minor extangential quadrilateral.   There is a sagelet to play with the values of $k$ and the $x$ coordinates of $B$ and $D$ at \textit{sagelets.org}.

\begin{center}
\includegraphics[scale=.3]{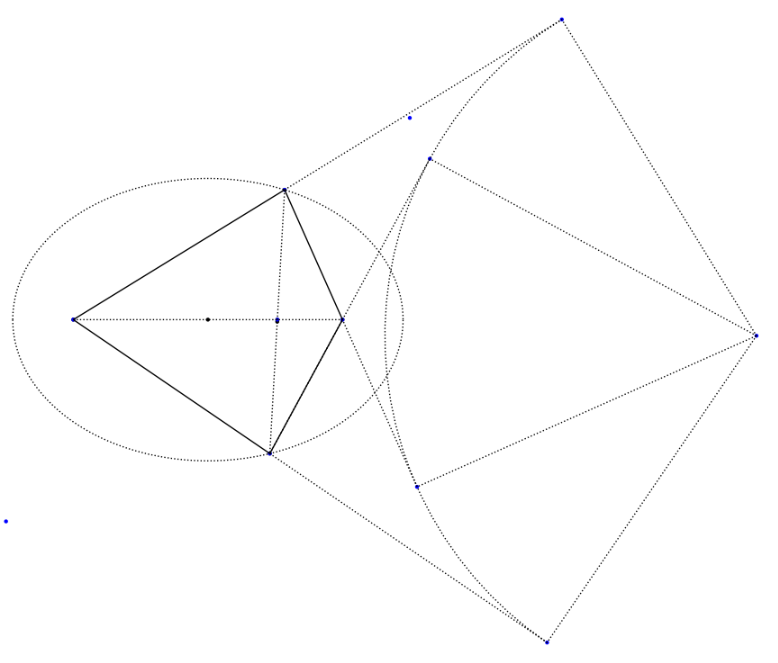}\\\textbf{The ellipse of a major extangential quadrilateral}
\end{center}

Each of the equations for the tangential and extangential quadrilaterals is a very ugly expression in the parameters $t_1,t_2,t_3, t_4$ equating the sum of two square roots with the sum of two other square roots. Since all four parameters are present in each equation, their graphs are not products, so we can't simply draw one cross section for each to see whether there is separation.  But we do know that their graphs are 3 dimensional and we can draw their 2-dimensional cross sections in one of the parameters.  Consider the  $t_2=.5$ cross sections of each equation drawn in the $t_1$, $t_3$, $t_4$ cube shown in the figures \textbf{View 1} and \textbf{View 2}.  The cross sections for the first equation lie in the top half of the cube and converge \textbf{down} to $t_2=1, t_4=0.5$ and the cross sections for equation 2 le in the bottom half and converge \textbf{up} to $t_2=1, t_4=.5$.

\begin{center}
\begin{tabular}{cc}
\includegraphics[scale=.4]{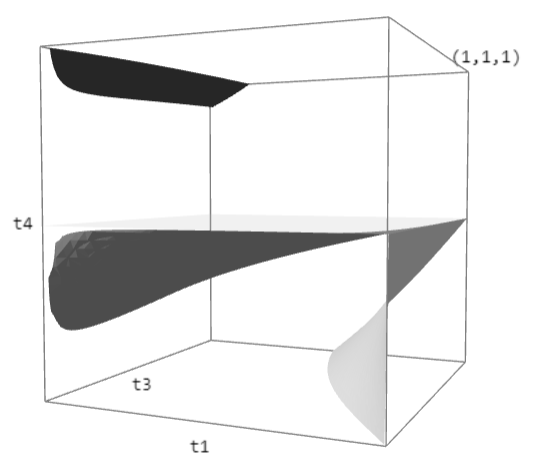} & \includegraphics[scale=.4]{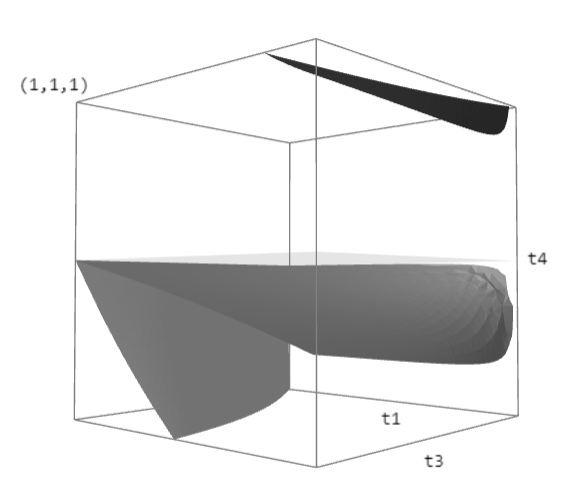}\\
\textbf{View 1} & \textbf{View 2}
\end{tabular}

\end{center}

Do the tangential (respectively extangential) quadrilaterals separate $\mathbb{Q}$?  A plausibility argument for separation analogous to the one for separation by cyclic quadrilaterals can be made here.
 For each $ABCD$ in $\mathbb{Q}$, let $\text{Tancir}(ABCD)$ be the unique circle which is tangent to the three sides $DA$,$AB$, and $BC$ of $ABCD$.  Decompose $\mathbb{Q}$
into three disjoint sets: $\mathbb{Tin}$,$\mathbb{Tan}$,$\mathbb{Tout}$ where $ABCD\in\mathbb{Tin}$,$\mathbb{Tan}$,or $\mathbb{Tout}$ according to whether the line $CD$ meets $\text{Tancir}(ABCD)$ in 0, 1 , or 2 points. So $\mathbb{Tan}$ is the class of cyclic quadrilaterals, and it is clear that if a quadrilateral from $\mathbb{Tin}$ is deformed along an arc of quadrilaterals to one in $\mathbb{Tout}$, then it must pass through $\mathbb{Tan}$.  The argument for extangential quadrilateral is analogous. The figures \textbf{View 1} and \textbf{View 2} support these plausibility arguments.  For any value of $t_2<1$, the graphs of the equations separate the the $t_2$ cross section of the $t_1,t_3,t_4$-cube  into three disjoint open sets.

\begin{center}
\begin{tabular}{ccc}
\includegraphics[scale=.2]{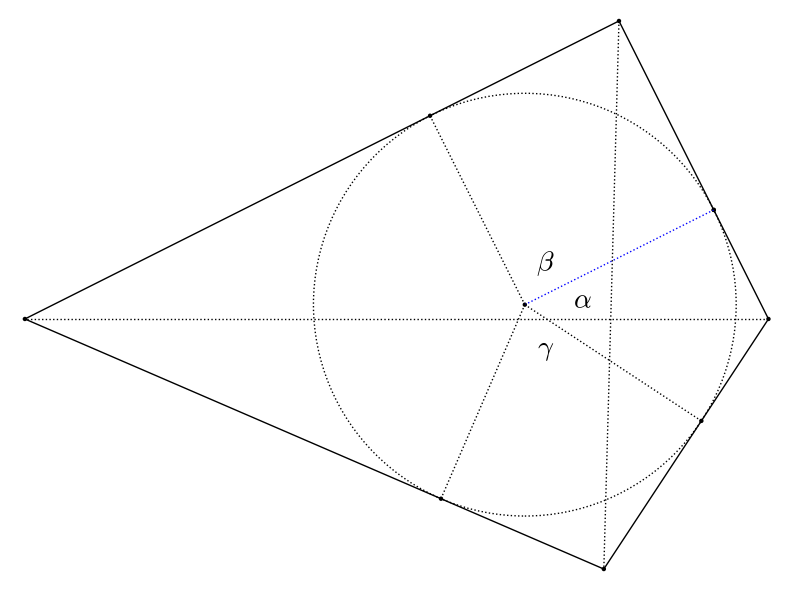} & \includegraphics[scale=.2]{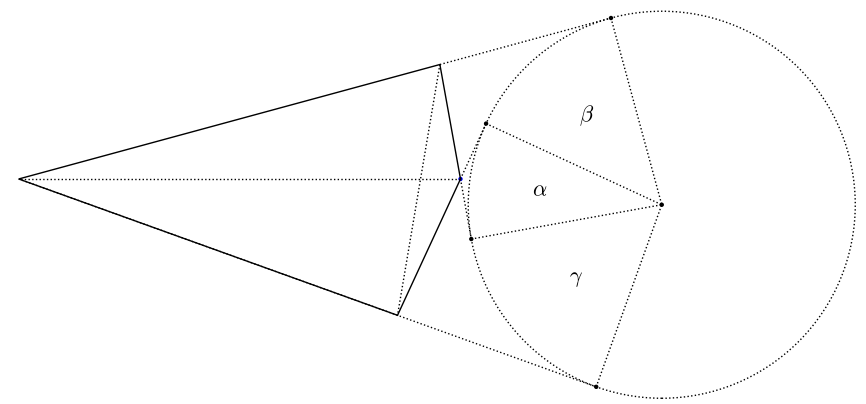} & \includegraphics[scale=.2]{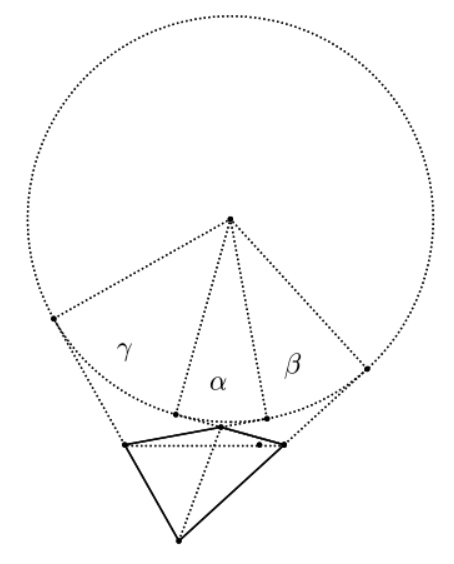}\\
\textbf{Tangential } & \textbf{Major Extangential }& \textbf{Minor Extangential  }
\end{tabular}
\end{center}

Since we can't solve (symbolically) either defining equation  for any one of the variables, it is difficult to use the equations to determine a specific tangential or ex-tangential quadrilateral.  One could resort to specifying 3 values and using numerical methods to determine the value of the 4th (if it exists).  This method produces unreliable results for us. However, there is another parametrization which is much less sensitive, involving the three angles $\alpha,\,\beta,\,\gamma$ shown in figures \textbf{Tangential}, \textbf{Major Extangential }, and \textbf{Minor Extangential}.

There is a sagelet implementing this parameterization at \textit{https://sagelets.cocalc.com/QuadSpace.html}. With it, you can supply $\alpha, \beta, \gamma$, and see it's picture and it's location in the graph.

\section{Questions}

In \cite{Eberhart}, we made the conjecture that \textit{if $P$ is a convex polygon with $2n+1$ vertices, then it has at most $2n+1$ quadrisections.}.   In efforts to prove this, we developed formulas for calculating the area of the upper right quadrant of each possible quadrisection. (A possile quadrisection is a pair of perpendicular lines each bisecting the area of the polygon.  There is one for each angle between 0 and $\pi/2$ radians.)   So far those efforts have failed to prove or disprove the conjecture.   We didn't make the same conjecture for quadrisections of polygons with an even number of sides because the square, having central symmetry has the property that each possible quadrisection is a quadrisection.  And we have found a tall isosceles trapezoid with exactly 5 quadrisections, and have included a picture of it and the graph of the area of the first quadrant minus one-fourth the area of the quadrilateral as the possible quadrisection rotates through 90 degrees.

\begin{center}
\begin{tabular}{cc}
\includegraphics[scale=.4]{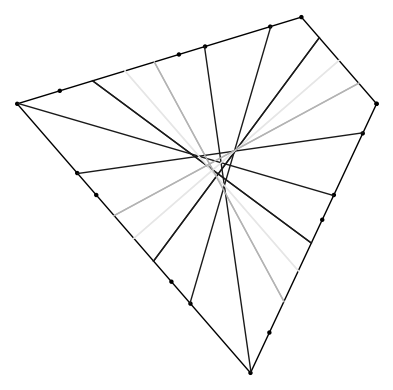} & \includegraphics[scale=.4]{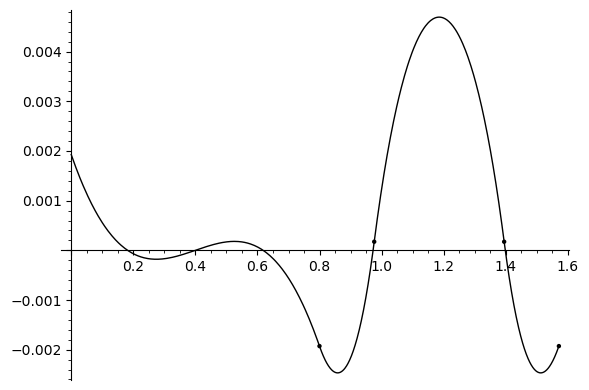}\\
\textbf{ Quardrilateral with 5 quadrisections. }& \textbf{ Graph of its area function. }
\end{tabular}

\end{center}

\textbf{Question:} \textit{Are there quadrilaterals with more than five but not infinitely many quadrisections?}  My guess is no, but that is only a guess, based on many failed searches for one with more than five.


\begin{bibdiv}
\begin{biblist}


\bib{Josefsson}{article}{
title={On the classification of convex quadrilaterals},
author={Josefsson, Martin},
date={March, 2016},
journal={The Mathematical Gazette, Vol 100, Issue 547, pp.68-85}}

\bib{Ahtziri}{article}{
title={Spaces of special quadrilaterals},
author={Gonzalez, Ahtziri and Lopez-Lopez, Jorge L.},
date={2019},
Journal={Bull. Aust. Math Soc. 100, pp. 155-167}
}

\bib{Eberhart}{article}{
title={Revisiting the quadrisection problem of Jacob Bernoulli.},
author={Eberhart, Carl},
date={Nov 21, 2016},
Journal={Forum Geometricorum, Vol 18, (2018) pp 7-16, found at \url{https://forumgeom.fau.edu/FG2018volume18/FG201802index.html }}
}



\end{biblist}
\end{bibdiv}
\end{document}